\documentclass[12pt]{amsart}
\usepackage{amssymb,amsmath, amsthm, enumerate, hyperref}
\usepackage[mathscr]{euscript}
\usepackage[a4paper]{geometry}
\usepackage[T1]{fontenc} 
\usepackage[utf8]{inputenc}   

\theoremstyle{plain}
\newtheorem{theorem}{Theorem}[section]
\theoremstyle{definition}

\newtheorem*{Acknowledgements}{Acknowledgement}
\numberwithin{equation}{section}

\def\R {\mathbb{R}}

\def\C {\mathbb{C}}

\def \au {\rm}
\def \ti {\it}
\def \jou {\rm}
\def \bk {\it}
\def \no#1#2#3 {{\bf #1} (#3), #2.}
\def \eds#1#2#3 {#1, #2, #3.}
\newcommand{\T}{(T(t))_{t\ge 0}}
\newcommand{\Ta}{(T_\alpha(t))_{t\ge 0}}
\DeclareMathOperator*{\re}{Re}

\begin{document}

\title[Gearhart--Pr\"uss theorem for bounded semigroups]
{A short elementary proof of the Gearhart--Pr\"uss\\
theorem for bounded semigroups}

\author[F. Dell'Oro]
{Filippo Dell'Oro}
\address{Politecnico di Milano -- Dipartimento di Matematica
\newline\indent
Via Bonardi 9, 20133 Milano, Italy}
\email{filippo.delloro@polimi.it}

\author[D. Seifert]
{David Seifert}
\address{Newcastle University -- School of Mathematics, Statistics and Physics
\newline\indent
Herschel Building, Newcastle upon Tyne, NE1 7RU, UK}
\email{david.seifert@ncl.ac.uk}

\subjclass[2020]{47D06, 47A10, 35B35}
\keywords{Gearhart--Pr\"uss theorem, bounded semigroup, exponential stability}

\begin{abstract}
We present a short elementary proof of the Gearhart--Pr\"uss 
theorem for bounded $C_0$-semigroups on Hilbert spaces. 
\end{abstract}

\maketitle

\section{Introduction}

\noindent Let $X$ be a complex Hilbert space and let $\T$ be a $C_0$-semigroup on $X$ with infinitesimal generator $A$. 
The semigroup $\T$ is said to be \emph{exponentially stable} if there exist $\omega>0$ and $M\geq1$ such that 
\begin{equation}
\label{exp-stab-def}
\|T(t)\|\le M {\rm e}^{-\omega t},\quad\, t\ge0.
\end{equation}
Exponential stability of $C_0$-semigroups is a classical notion, but
it remains one of the most important types of stability and it arises 
in numerous applications. Some of these
may be found in the references given below, although we make no attempt in this note 
at giving even a partial overview of the relevant literature. 
Since one typically has no explicit knowledge of the operators $T(t)$ for $t>0$, it is desirable to have a 
criterion for exponential stability purely in terms of the (spectral) properties of the infinitesimal generator $A$.
We denote, using standard notation,  by $\rho(A)$ the resolvent set of $A$ and by $R(\lambda,A)=(\lambda- A)^{-1}$ the 
resolvent operator at $\lambda\in\rho(A)$. We also denote by $i\R$ the imaginary axis in the complex plane. 
The following result is a version of the celebrated Gearhart--Prüss theorem
for bounded $C_0$-semigroups. Recall that a $C_0$-semigroup $\T$ is said to be \emph{bounded} if
\begin{equation}\label{bound-def}
\sup_{t\ge0}\|T(t)\|<\infty,
\end{equation}
and in this case the open right half-plane $\C_{+} = \{ \lambda \in \C : \re \lambda > 0\}$
is necessarily contained in $\rho(A)$;\ see e.g.\ \cite[Theorem~II.1.10]{ENG}.

\begin{theorem}
\label{thm:GP}
Let $X$ be a complex Hilbert space and let $\T$ be a bounded $C_0$-semigroup on $X$ 
with infinitesimal generator $A$. 
Then $\T$ is exponentially stable if and only if $i\R\subseteq\rho(A)$
and $\sup_{s\in\R}\|R(is,A)\|<\infty$.
\end{theorem}

Note that, by continuity of the map $\lambda\mapsto \|R(\lambda,A)\|$  on $\rho(A)$,  
uniform boundedness of the resolvent operator on the imaginary axis is equivalent to the condition that 
$$\limsup_{|s|\to\infty}\|R(is,A)\|<\infty,$$
which is  commonly found in some parts of the literature.
Theorem \ref{thm:GP} is a straightforward consequence of a more general version of the {Gearhart--Pr\"uss theorem},
which states that a (not necessarily bounded) $C_0$-semigroup $\T$ on a 
Hilbert space is exponentially stable if and only if  $\C_{+}  \subseteq\rho(A)$ and 
$\sup_{\lambda \in \C_{+}}\|R(\lambda, A)\| <\infty.$
The latter result originated in the work of Gearhart~\cite{GER} and the subsequent 
contributions (in chronological order) by Monauni~\cite{MON}, Herbst~\cite{HER}, 
Howland~\cite{HOW}, Pr\"uss~\cite{PRU}, Huang~\cite{HUA} and Greiner~\cite{ARE}.
More detailed discussions and historical accounts may be found for instance in \cite{BATTYBOOK,CHI,CST,CRA,ENG}.
One of the most common proofs of the general version of the Gearhart--Pr\"uss theorem, given 
for instance in the well-known textbook \cite[Chapter~V]{ENG},
requires the use of several non-trivial (albeit standard) tools, 
such as a special form of the complex inversion formula for the Laplace transform and Cauchy's integral theorem. 
Elegant short proofs based on the Datko--Pazy theorem may be found in \cite[Chapter~5]{BATTYBOOK} 
and \cite{WEISS}, while other short proofs based on the Paley--Wiener 
theorem are given in \cite[Chapter 3]{VAN} and \cite[Chapter~5]{CZ}. 
A proof relying on the theory of  so-called evolution semigroups may be found in \cite{LAT}.
We refer to the more recent \cite[Theorem~4.7]{BCT} for a different short and elementary proof of the 
Gearhart--Prüss theorem specifically for bounded semigroups (this is not stated explicitly 
by the authors of \cite{BCT} but follows at once on taking $B$, in their notation, to be the identity operator).
We further mention the article of Wyler \cite{WYL}, which contains the first explicit statement (to our knowledge)
of the sufficiency part of Theorem \ref{thm:GP}. 
Other relatively recent works relating to the Gearhart--Prüss theorem include \cite{HEL, WEI}.

The aim of this note is to give another short and elementary proof of the Gearhart--Pr\"uss 
theorem for bounded semigroups, as formulated in Theorem~\ref{thm:GP}. 
Our proof, to be given in Section~\ref{sec:proof} below, relies almost exclusively on a small number of
basic facts from elementary semigroup theory, such as the semigroup property, 
the resolvent identity and the integral representation of the resolvent of a semigroup generator. 
In particular, we require neither the Datko--Pazy theorem nor the Paley--Wiener theorem and,
in addition, we avoid using the adjoint semigroup, 
the inversion formula and even Cauchy's integral theorem. 
Arguably the only non-elementary result used in our proof is Plancherel's theorem 
for square-integrable functions taking values in a Hilbert space; 
see e.g.\ \cite[Appendix C]{ENG}. The latter,
in one way or another, plays a key role in all known proofs of the Gearhart--Prüss theorem. 
We conclude our paper in Section~\ref{sec:rem} with a discussion in the form of some remarks.

Our main motivation for writing this note is that many (if not most) modern applications 
of the Gearhart--Pr\"uss theorem to PDEs arise in the setting of 
bounded or indeed contractive semigroups (that is, $C_0$-semigroups for which the supremum 
in \eqref{bound-def} equals one). As a result, a short elementary 
proof of Theorem~\ref{thm:GP} may be of interest to certain parts of both the PDE and  the 
semigroup communities.
We moreover hope that our  proof will be of value to those giving introductory 
lecture courses on the theory of $C_0$-semigroups. 

\section{Proof of Theorem \ref{thm:GP}}
\label{sec:proof}

\noindent
Since there will be no risk of confusion, 
throughout this section we shall simply write $R(\lambda)$ for $R(\lambda,A)$ when $\lambda\in\rho(A)$.

The necessity part of Theorem \ref{thm:GP} is standard. Indeed, if $\T$ is exponentially stable 
then we may find  $\omega>0$ and $M\ge1$ such that \eqref{exp-stab-def} holds. 
Moreover, for $\lambda\in\C$ with $\re\lambda>-\omega$, we have $\lambda\in\rho(A)$ and 
 the resolvent  admits the integral representation
$$
R(\lambda)x=\int_0^\infty {\rm e}^{-\lambda t} T(t)x\,\mathrm{d} t,\quad\, x\in X;
$$
see e.g.\ \cite[Theorem~II.1.10]{ENG}.
In particular, $i\R\subseteq\rho(A)$ and a straightforward estimate 
yields $\sup_{s\in\R}\|R(is)\|\le M\omega^{-1}<\infty$.

In order to prove the sufficiency part of Theorem \ref{thm:GP} we note that,
as a consequence of the semigroup property, 
it suffices to show that $\lim_{t\to\infty}\|T(t)\|=0$; 
see e.g.\ \cite[Proposition~V.1.7]{ENG}.
Suppose then that $i\R\subseteq\rho(A)$ and $\sup_{s\in\R}\|R(is)\|<\infty$, and let $K$ 
be the finite supremum in \eqref{bound-def}.  
Given $\alpha>0$, we consider the rescaled semigroup $\Ta$ 
given by $T_\alpha(t) = e^{-\alpha t} T(t)$ for $t\ge0$.
For every $x\in X$ and $t>0$ we have
\begin{equation}\label{eq:Ta_bd}
\|T_\alpha(t)x\|^2= \frac1t\int_0^t \|T_\alpha(t-\tau) T_\alpha(\tau) x\|^2 \, \mathrm{d}\tau 
\leq \frac{K^2}{t} \int_0^{\infty}\|T_\alpha(\tau)x\|^2 \, \mathrm{d}\tau.
\end{equation}
Now let $C=\sup_{s\in\R}\|R(is)\|$  and $\omega=C^{-1}$.
It follows from the resolvent identity that, for $\alpha>0$ and $s\in\R$, 
$$\|R(\alpha+i s)\| =\|R(i s )-\alpha R(\alpha+i s) R(i s) \| \leq C+ C\alpha  \|R(\alpha+i s) \|,$$
and hence $\|R(\alpha+i s)\|\le (\omega- \alpha)^{-1}$ for all $\alpha\in(0,\omega)$ and $s\in\R$.
By another application of the resolvent identity
we obtain
\begin{equation*}\label{eq:res_bd}
\|R(\alpha +i s)x\| = \|R(\omega +i s)x+(\omega-\alpha) R(\alpha +i s) R(\omega +i s)x\|
\leq 2 \|R(\omega+ i s)x\|
\end{equation*}
for all $\alpha\in(0,\omega)$, $s\in\R$ and $x\in X$.
If  $\alpha>0$ and $x\in X$, and if we extend the semigroup $\Ta$ by zero to $\R$, 
then by the integral representation of the resolvent the function $s\mapsto R(\alpha + i s)x$ 
is the Fourier transform of the function
$t\mapsto T_\alpha(t) x$. Thus, applying Plancherel's theorem (twice) we deduce that, for $\alpha\in(0,\omega)$,
\begin{equation}\label{eq:int_bd}
\begin{aligned}
\int_0^{\infty} \|T_\alpha(\tau)x\|^2\, \mathrm{d}\tau 
&= \frac{1}{2\pi}\int_{-\infty}^{\infty} \|R(\alpha + i s)x \|^2\, \mathrm{d}s
\leq \frac{2}{\pi} \int_{-\infty}^{\infty}\|R(\omega + i s)x\|^2\, \mathrm{d}s\\
&= 4 \int_0^{\infty} \|T_{\omega}(\tau)x\|^2\, \mathrm{d}\tau\le c^2\|x\|^2,
\end{aligned}
\end{equation}
 where $c=K(2C)^{1/2}$. Combining \eqref{eq:Ta_bd} and \eqref{eq:int_bd} 
gives $\|T_\alpha(t)\| \leq Kc\,t^{-1/2}$ for all $\alpha\in(0,\omega)$ and $t>0$, and letting $\alpha\to0^{+}$ we get
$$\|T(t)\|\leq \frac{Kc }{ t^{1/2}}\to0, \quad\, t\to\infty.$$ 
Thus $\T$ is exponentially stable, as required.
\qed

\section{Concluding Remarks}
\label{sec:rem}

\noindent (a) In contrast to the arguments given in \cite{BCT, ENG}, the above proof of Theorem~\ref{thm:GP} does not
make use of the adjoint semigroup $(T(t)^*)_{t\ge0}$. 
However, one may alternatively observe that the 
same arguments that led to \eqref{eq:int_bd} when applied to the adjoint semigroup yield 
\begin{equation}\label{eq:int_bd2}
\int_0^{\infty} \|T_\alpha^*(\tau)y\|^2\, \mathrm{d}\tau \le c^2\|y\|^2
\end{equation}
for all $\alpha\in(0,\omega)$ and $y\in X$, where $c$ is 
as above and $T_\alpha^*(t)=e^{-\alpha t}T(t)^*$ for all $t\ge0$. 
Now let $x,y\in X$, $t>0$ and $\alpha\in(0,\omega)$. 
Combining \eqref{eq:int_bd} and \eqref{eq:int_bd2} with an application of the Cauchy--Schwarz inequality, we obtain
$$|\langle T_\alpha (t)x,y\rangle|
=\frac1t\int_0^t|\langle T_\alpha (\tau)x,T_\alpha^*(t-\tau)y\rangle|\,\mathrm{d}\tau\le \frac{c^2}{t}\|x\|\|y\|,$$
and hence $\|T_\alpha(t)\|\le c^2\,t^{-1}$. 
Letting $\alpha\to0^{+}$ gives
$\|T(t)\|\leq c^2 \, t^{-1}$ for all $t>0$,
and the result once again follows. Note that this slightly more involved approach produces the (better) decay
rate $\|T(t)\|={\rm O}(t^{-1})$ as $t\to\infty$.\smallskip

\noindent(b)
It follows from the Datko--Pazy theorem \cite[Theorem~5.1.2]{BATTYBOOK} 
that a (not necessarily bounded) $C_0$-semigroup $\T$ with infinitesimal generator $A$ 
acting on a complex \emph{Banach space} $X$ is exponentially stable 
if and only if there exists $p\in[1,\infty)$ such that $T*f\in L^p(0,\infty;X)$ for all $f\in L^p(0,\infty;X)$,
where
$$
(T*f)(t) = \int_0^t T(t-s) f(s) \,\mathrm{d}s,\quad \, t\ge0.
$$
If this is the case, an application of the closed graph theorem shows 
that the map $f\mapsto T*f$ is automatically bounded on $L^p(0,\infty;X)$. 
An equivalent condition is that the operator-valued 
map $m$ sending $\lambda\in\C_+$ to $ R(\lambda,A)$ is a 
well-defined $L^p(0,\infty;X)$-{Laplace multiplier}; see \cite[Remark~4.8]{BCT}.
In particular, the $C_0$-semigroup $\T$ is exponentially stable 
if and only if $m$ is an $L^p(0,\infty;X)$-{Laplace multiplier} for some $p\in[1,\infty)$, and this result may be viewed
as a generalisation of the Gearhart--Prüss theorem to Banach spaces.
If $X$ is a Hilbert space, then it follows from Plancherel's theorem 
that $m$ is an $L^2(0,\infty;X)$-Laplace multiplier whenever it is well defined 
and bounded, thus recovering the classical Gearhart--Prüss theorem. 

As noted in \cite[Remark~4.8]{BCT}, when $\T$ is bounded the argument used to prove 
\cite[Theorem~4.7]{BCT} can be adapted to give an alternative 
proof of the aforementioned version of the Gearhart--Prüss theorem on Banach spaces. \smallskip

\noindent (c)
It is worth remarking that the Gearhart--Prüss theorem does not extend to arbitrary Banach spaces, 
or even to $L^p$-spaces with $p\in(1,\infty)$ and $p\neq 2$. 
We refer to \cite[Example~5.2.2]{BATTYBOOK} and \cite[Comments V.1.12]{ENG} for suitable examples. However, a version of the Gearhart-Prüss theorem on $K$-convex spaces may be found in~\cite{Arn20}.

\begin{Acknowledgements}
The authors thank Vittorino Pata  for useful discussions.
\end{Acknowledgements}



\begin{thebibliography}{00}

\bibitem{BATTYBOOK}
{\au W. Arendt, C.J.K. Batty, M. Hieber and F. Neubrander},
{\bk Vector-valued Laplace Transforms and Cauchy Problems},
\eds{Birkh\"auser}{Basel}{second edition, 2011}

\bibitem{ARE}
{\au W. Arendt, A. Grabosch, G. Greiner, 
U. Groh, H.P. Lotz, U. Moustakas, R. Nagel,
F. Neubrander and U. Schlotterbeck}, 
{\ti One-parameter Semigroups of Positive Operators}, 
\eds{Lecture Notes in Mathematics, 1184. Springer-Verlag}{Berlin}{1986}

\bibitem{Arn20}
{\au L. Arnold},
{\ti $\gamma$-boundedness of $C_0$-semigroups and their $H^\infty$-functional calculi},
{\jou Studia Math.}
\no{254}{77--108}{2020}

\bibitem{BCT}
{\au C.J.K. Batty, R. Chill and Yu. Tomilov},
{\ti Fine scales of decay of operator semigroups},
{\jou J.\ Eur.\ Math.\ Soc.\ (JEMS)}
\no{18}{853--929}{2016}

\bibitem{CHI}
{\au C. Chicone and Yu. Latushkin},
{\bk Evolution Semigroups in Dynamical Systems and Differential Equations},
\eds{Amer.\ Math.\ Soc.}{Providence}{1999}

\bibitem{CST}
{\au R. Chill, D. Seifert and Yu. Tomilov}, 
{\ti Semi-uniform stability of operator semigroups and energy decay of damped waves}, 
{\jou Philos.\ Trans.\ Roy.\ Soc.\ A}
\no{378}{20190614, 24 pp}{2020}

\bibitem{CRA}
{\au D. Cramer and Yu. Latushkin},
{\bk Gearhart--Pr\"uss Theorem in Stability for Wave Equations:\ a Survey},
\eds{Lecture Notes in Pure and Appl.\ Math., 234}{Dekker, New York}{2003} 

\bibitem{CZ}
{\au R.F. Curtain and H.J. Zwart},
{\bk An Introduction to Infinite-dimensional Linear System Theory},
\eds{Springer-Verlag}{New York}{1995}

\bibitem{ENG}
{\au K.-J. Engel and R. Nagel},
{\bk One-parameter Semigroups for Linear Evolution Equations},
\eds{Springer-Verlag}{New York}{2000}

\bibitem{GER}
{\au L. Gearhart},
{\ti Spectral theory for contraction semigroups on Hilbert space},
{\jou Trans.\ Amer.\ Math.\ Soc.}
\no{236}{385--394}{1978}

\bibitem{HEL}
{\au B. Helffer and J. Sjöstrand},
{\ti From resolvent bounds to semigroup bounds},
Actes du colloque d'Evian, 2009. https://arxiv.org/abs/1001.4171v1.

\bibitem{HER}
{\au W. Herbst},
{\ti The spectrum of Hilbert space semigroups}, 
{\jou J.\ Operator Theory} 
\no{10}{87--94}{1983}

\bibitem{HOW}
{\au J.S. Howland}, 
{\ti On a theorem of Gearhart}, 
{\jou  Integral Equ. Oper. Theory} 
\no{7}{138--142}{1984} 

\bibitem{HUA}
{\au F.L. Huang},
{\ti Characteristic conditions for exponential stability 
of linear dynamical systems in Hilbert spaces},
{\jou Ann.\ Differential Equations}
\no{1}{43--56}{1985}

\bibitem{LAT}
{\au Yu. Latushkin and S. Montgomery-Smith},
{\ti Evolutionary semigroups and Lyapunov theorems in Banach spaces},
{\jou J.\ Funct.\ Anal} 
\no{127}{173--197}{1995} 

\bibitem{MON}
{\au L.A. Monauni},
{\ti On the abstract Cauchy problem and the generation problem 
for semigroups of bounded operators},
\eds{Control Theory Centre Report No.\ 90}{Warwick}{1980}

\bibitem{VAN}
{\au J. van Neerven},
{\bk The Asymptotic Behaviour of Semigroups of Linear Operators}, 
\eds{Operator Theory:\ Advances and Applications, 
88. Birkhäuser Verlag}{Basel}{1996}

\bibitem{PRU}
{\au J. Pr\"{u}ss},
{\ti On the spectrum of $\text{C}_0$-semigroups},
{\jou Trans.\ Amer.\ Math.\ Soc.}
\no{284}{847--857}{1984}


\bibitem{WEI}
{\au D. Wei},
{\ti Diffusion and mixing in fluid flow via the resolvent estimate},
{\jou Sci.\ China Math.}
\no{64}{507--518}{2021} 

\bibitem{WEISS}
{\au  G. Weiss}, 
{\ti Weak $L^p$-stability of a linear semigroup on a Hilbert space implies exponential stability}, 
{\jou J. Differential Equations} 
\no{76}{269--285}{1988}

\bibitem{WYL}
{\au A. Wyler},
{\ti Stability of wave equations with dissipative boundary conditions in a bounded domain}, 
{\jou Differential Integral Equations} 
\no{7}{345--366}{1994} 

\end{thebibliography}
\end{document}